\documentclass[12pt]{article}

\usepackage[a4paper,margin=1in,footskip=0.25in]{geometry}
\usepackage{color,soul}

\usepackage{latexsym,amssymb,lastpage}
\usepackage{graphicx,amsfonts}
\usepackage{times,mathptmx,bm,amsmath,amsthm}
\usepackage{dcolumn}

\newenvironment{pf}[1][Proof]{\noindent\textit{#1. } }{\hfill$\square$}
\renewcommand{\Re}{\mathbb{R}}

\newtheorem{prop}{Proposition}
\newtheorem{thm}{Theorem}
\newcommand{\diam}{\delta_{\footnotesize{\mbox{diam}}}}
\renewcommand{\theequation}{\thesection.\arabic{equation}}
\numberwithin{equation}{section}

\newcommand{\sH}{\mathcal{H}}
\newcommand{\sE}{\mathcal{E}}

\begin{document}
\title{Open Problem Statement: \\
Minimal Distortion Embeddings of Diversities in $\ell_1$}
\author{David Bryant and Paul Tupper}
\date{\today}
\maketitle

\vspace{-.5cm}

\section{Statement of Problem}

A diversity is a pair $(X,\delta)$ where $X$ is a set and  $\delta$ is a function taking finite subsets of $X$ to $\mathbb{R}$ and satisfying the following axioms: \cite{bryant2012} \\
(i) $\delta(A) \geq 0$ for all $A \subset X$ \\
(ii) $\delta(A)= 0$ if and only if $|A|\leq 1$ \\
(iii) $\delta(A \cup B) \leq \delta(A \cup C) + \delta(B \cup C)$ whenever $C \neq \emptyset$. \\
Diversities are an extension of the concept of a metric space 
since for any diversity if we define $d(x,y) = \delta(\{x,y\})$ for $x, y\in X$ then $(X,d)$ is a metric space \cite{bryant2012}. Work on the mathematics of diversities has generated new theoretical results, applications, and unexpected links \cite{deza2016generalizations,EspinolaETAL2014a,HerrmannETAL2012a,Poelstra2013b,Piatek2014a,KumarETAL2015a,KirkETAL2014a,Steel2014a}.

Just as there is an $\ell_1$ metric on $\mathbb{R}^k$ there is an $\ell_1$ diversity on $\mathbb{R}^k$ defined as
\[
\hat{\delta}(A) := \sum_{i=1}^k \max_{a,b \in A} |a_i - b_i|
\]
for all finite subsets $A$ of $\mathbb{R}^k$.

Let $(X_1,\delta_1)$ and $(X_2,\delta_2)$ be two diversities.   Let $\phi$ be a map from $X_1$ to $X_2$. We say $\phi$ has distortion $c$ if there are are $c_1, c_2 >0$ such that $c=c_1 c_2$ and
\[
\frac{1}{c_1} \delta_1(A) \leq \delta_2 (\phi(A)) \leq c_2 \delta_1(A)
\]
for all finite $A \subseteq X_1$ \cite{bryant2014}. 


{\bf Problem:} Find a bound $f(n)$ such that for any  diversity $(X,\delta)$ with $|X|=n$ there is a map $\phi$ from $(X,\delta)$ into an $\ell_1$ diversity $(\mathbb{R}^d,\hat{\delta})$  with distortion at most $f(n)$. \pagebreak


As we explain below, the best general bound we have so far is $\mathcal{O}(n)$, though this is likely far from optimal. We know that the bound cannot be better than $\Omega(\log n)$.


Our problem is the diversity analogue of a well-known problem in the theory of metric spaces: Find the minimal distortion embedding of a $n$-point metric space into $(\mathbb{R}^k,\|\cdot\|_1)$ for any $k$. Due to a result of Bourgain \cite{Bourgain85}, it is known that the distortion bound in the metric spaces case is $\mathcal{O}(\log n)$ and that this can be attained with dimension $k=\mathcal{O}(\log n)$. This bound is tight \cite{Linial95}.

We know that for the diversity case it is not possible to embed a general $n$-point diversities with distortion $\mathcal{O}(\log n)$ into $\ell_1$ with dimension $k=\mathcal{O}(\log n)$ \cite{bryant2014}. However, once we remove the restriction on the dimension, the $\mathcal{O}(\log n)$ distortion has not been ruled out yet.

The following  argument, due to P. Wu, gives an upper bound for the worst distortion of $\mathcal{O}(n)$.   
\begin{thm} \label{PetersEmbedding}
For all diversities $(X,\delta)$ where $|X|=n$ there is an embedding $\phi \colon X \rightarrow \ell_1^n$ with distortion at most $n$.
\end{thm}
\begin{pf}
Let $(X,d)$ be the induced metric of $(X,\delta)$. 
Letting $X = \{x_1, \ldots, x_n\}$.  Define $\phi$ by
\[
\phi(x) = ( d(x,x_1),d(x,x_2),\ldots,d(x,x_n) ).
\]
Then for all $A \subseteq X$ we have
\[
n \delta(A) \geq \sum_{i=1}^n \max_{a,b \in A} d(a,b) \geq \sum_{i} \max_{a,b \in A} d(a,x_i) - d(b,x_i) = \hat{\delta}(\phi(A)),
\]
where the first inequality follows from the monotonicity of $\delta$ and the second follows from the triangle inequality for diversities.
Also, for all $A \subseteq X$ we have
\[
\hat{\delta}(\phi(A)) \geq \sum_{x \in A}  \max_{a,b \in A} d(a,x) - d(b,x) = \sum_{x \in A} \max_{a \in A} d(a,x) \geq \sum_{x \in A} d(a_0,x) \geq \delta(A),
\]
where $a_0$ is any element of $A$. These two inequalities give distortion at most $n$.
\end{pf}

There are metric spaces that require $\Omega(\log n)$ distortion \cite{Linial95}. If we take a diversity that has this metric as its induced metric, then we cannot embed that diversity with less distortion.

In Section~\ref{sec:specCase} we give a survey of special classes of diversities for which we have upper bounds better than $\mathcal{O}(n)$. In Section~\ref{sec:failAttempts} we show how two of the techniques that were used for the metric embedding problem fail for diversities. Section~\ref{sec:otherIdeas} lists some ideas for possible directions to follow.

\section{Special Cases}\label{sec:specCase}

There are a few special cases where we have tighter upper bounds on the worst case distortion.
See \cite{bryant2014} for examples of diversities which can be embedded into $\ell_1$ with no distortion, as well as characterisations of these diversities.

\subsection{Diameter Diversities}

Every metric space $(X,d)$ induces a diversity $(X,\diam)$, where $\diam$ is the diameter diversity defined by
\[
\diam(A) = \max_{a,b \in A} d(a,b).
\]
$(X,\diam)$ is a diversity, and is in fact the minimal diversity that has $d$ as its induced metric \cite{bryant2012}

Let $\ell^k_1$ denote $\mathbb{R}^k$ equipped with the $\ell_1$ diversity $\hat{\delta}$.
In \cite[Prop. 12]{bryant2014} we show that any diameter diversity on $n$ points can be embedded into $\ell_1^k$ with distortion $\mathcal{O}(\log^2 n)$ where $k = \mathcal{O}(\log n)$. The result uses the $\ell_1$ embedding result of \cite{Linial95} along with the fact that there is at most a ratio of $k$ between the diameter diversity and the $\ell_1$ diversity on $\ell_1^k$.

One important example of a diameter diversity is the {\em discrete diversity}, given by 
\[\delta_\rho(A) = \begin{cases} 1 & |A|>1 \\ 0 & \mbox{ otherwise.} \end{cases}\]

\subsection{Steiner Diversities}

Consider any metric space $(X,d)$.  View the metric space as a complete weighted graph with vertices $X$ and the weight of edge $(u,v)$ to be $d(u,v)$.   The Steiner diversity $\delta_S(A)$ of a finite set $A$  is the minimal total weight of a  tree in the graph whose leaves include the points in $A$.  $(X,\delta_S)$ is known as the Steiner diversity of $(X,d)$, and is the maximal diversity that has $(X,d)$ as its induced metric \cite{bryant2014}.

A simple argument based on a metric approximation result due to \cite{fakchar04} gives the following.
We note that the result of \cite{fakchar04} was a refinement of an earlier  result of Bartal \cite{bartal96,bartal98}.
\begin{thm} \label{SteinerBound}
Let $(X,d)$ be  a metric space with $|X|=n$ and let $(X,\delta_S)$ be its Steiner diversity. Then $(X,\delta_S)$ has an $\ell_1$ embedding with distortion $\mathcal{O}(\log n)$.
\end{thm}
\begin{pf}
By \cite{fakchar04} there is a  bijection $\phi$ from $X$ to the vertices of a random tree $\tau$ with vertex set $X$ such that for all $u,v \in X$ we have $d(u,v) \leq d_\tau(u,v)$ and $\mathbb{E}[d_\tau(u,v)] \leq (\log n)\, d(u,v)$. (Here $\mathbb{E}$ denotes expectation.)

For all  $A \subseteq X$, there is a tree $T = (X_{\phi(A)},E_{\phi(A)})$ such that $T$ is a subtree connecting $\phi(A)$ in $\tau$ and $\sum_{\{u,v\} \in E_{\phi(A)}} d_\tau(u,v) = \delta_\tau(A)$. Then 
\[\delta_S(A) \leq \sum_{\{u,v\} \in E_{\phi(A)}} d(u,v) \leq  \sum_{\{u,v\} \in E_{\phi(A)}} d_\tau(u,v) = \delta_\tau(A) .\]

Likewise, for all $A \subseteq X$, there is a tree $T =(V_A,E_A)$ such that 
\[\delta_S(A) = \sum_{\{u,v\} \in E_A} d(u,v) \]
(by the definition of the Steiner diversity)
which gives 
\begin{align*}
\mathbb{E}[\delta_\tau(A)] & \leq \mathbb{E}\left[ \sum_{\{u,v\} \in E_A} d_\tau(u,v) \right] \\
& \leq \mathcal{O}(\log n) \sum_{\{u,v\} \in E_A} d(u,v) \\
& = \mathcal{O}(\log n) \delta(A).
\end{align*}
Hence $\mathbb{E}[\delta_\tau(A)] \leq \mathcal{O}(\log n) \delta(A)$.

For each $\tau$ we have that $\delta_\tau$ can be embedded in $\ell_1$, see \cite{bryant2012}. Hence so can the expectation,
\[\hat{\delta} = \mathbb{E} \delta_{\tau} \]
and from the above we have that
\[\delta(A) \leq \hat{\delta}(A) \leq \mathcal{O}(\log n) \delta(A)\]
for all $A$. 
\end{pf}

\subsection{Hypergraph Steiner}

Let $\sH =(V,\sE,w)$ be a hypergraph with a weight function $w:\sE \rightarrow \Re$. The hypergraph Steiner diversity $\delta_H$ is defined as the minimum of
\[w(\sH') = \sum_{U \in \sE(\sH')} w(U)\]
over all connected sub-hypergraphs of $\sH$ with $A \subseteq V(\sH')$, with singletons defined to have diversity zero.

{\em Every} diversity on a finite set equals $\delta_H$ for some weighted hypergraph---just use $\sE = P(V)$ and $w(U) = \delta(U)$. We bound the embedding distortion for $\delta_H$ in terms of the sizes of hyperedges in $\sE(\sH)$.

\begin{thm}
Let $\delta_H$ be the Hypergraph Steiner diversity for $\sH = (V,\sE,w)$ and suppose that $k = \max\{|U|:U \in \sE(\sH)\}$. Then $\delta_H$ can be embedded in $\ell_1$ with $\mathcal{O}(k \log n)$ distortion.
\end{thm}
\begin{pf}
For each $U \in \sE$ let $T_U$ be a spanning tree of the connected graph with vertex set $U$, where each edge in $T_U$ has weight $w(U)$. Construct the graph $G = (V,E)$ where $E$ is the disjoint union
\[E = \biguplus_{U \in \sE}  E(T_U).\] 

Every connected subgraph of $G$ gives a connected sub-hypergraph of $\sH$ with lesser weight, and every connected sub-hypergraph of $\sH$ gives a  connected subgraph of $G$ with weight at most $(k-1)$ times larger. It follows that 
\[\delta_H(A) \leq \delta_S(A) \leq (k-1) \delta_H(A)\]
and the result now follows from Theorem~\ref{SteinerBound}.
\end{pf}

\subsection{Ball Diversities}

Let $(X,d)$ be any metric space. Define the diversity $(X,\delta_B)$ by letting $\delta_B(A)$ equal the diameter of the smallest ball containing $A$. We call this the {\em ball diversity}.

For any $A \subseteq X$ we have 
\[\diam(A) \leq \delta_B(A) \leq 2 \min_{a \in A} \max_{b \in A} d(a,b) \leq 2 \diam(A)\]
and so a bound of $\mathcal{O}(\log^2 n)$ distortion for embedding ball diversities in $\ell_1$ follows from the diameter diversity case.

\subsection{TSP Diversity}

Let $(X,d)$ be any metric space. Define the diversity $(X,\delta_{TSP})$ where $\delta_{TSP}(A)$ is the length of the shortest tour through points in $A$. We then have that 
\[\delta_{S}(A) \leq \delta_{TSP}(A) \leq 2 \delta_{S}(A)\]
so $\delta_{TSP}$ can be embedded with $\mathcal{O}(\log n)$ distortion.

\subsection{Partition diversities}

Given a partition $\pi$ of $X$ we define $\delta_\pi$ by
\[\delta_{\pi}(A) = \begin{cases} 1 & \mbox{ $A$ intersects at least two blocks of $\pi$; } \\ 0 & \mbox{ otherwise. } \end{cases}
\]
This reduces to a {\em split diversity} or {\em cut diversity} when $\pi$ has two blocks \cite{bryant2014}. In general, partition diversities are not $\ell_1$ embeddable, however they are examples of diameter diversities. Consequently, any non-negative linear combination of partition diversities can be embedded into $\ell_1$ with $\mathcal{O}(\log^2 n)$ distortion.

\subsection{Symmetric Diversities}

We say a diversity is symmetric if $\delta(A)=\delta(B)$ whenever $|A|=|B|$.   Perhaps surprisingly, symmetric diversities can be embedded into $\ell_1$ with a constant factor distortion \cite{bryant2016}. A consequence of this result is that if $\frac{\delta(A)}{\delta(B)} \leq K$ whenever $|A| = |B| > 1$ then $\delta$ can be embedded in $\ell_1$ with distortion 
$\mathcal{O}(K)$.


\section{Failed Attempts}\label{sec:failAttempts}

We have tried two different approaches to the problem, one based on the techniques used in Bourgain's proof, and the other based on Bartal's techniques. Both were unsuccessful, but we will explain where the problem lies in each case.

\subsection{Bourgain's Method}

To summarize Bourgain's method at a high level \cite{Bourgain85,Deza97}, suppose we are given a metric $(X,d)$ with $|X|=n$ and we want to embed it in $\ell_1^k$ for some $k$. We select a collection of maps $\phi_{d,A} \colon X \rightarrow \mathbb{R}$ for $A \subseteq X$, and a collection of weights $c_{d,A} \geq 0$. Then the $\ell_1$ metric is defined as
\[
\hat{d}(x,y)= \sum_{A \subseteq X} c_{d,A} | \phi_{d,A}(x) - \phi_{d,A}(y) |.
\]
The space $(X,\hat{d})$ can be embedded in $\ell_1^k$ where $k$ is the number of non-zero $c_{d,A}$. The $\phi_{d,A}$ in the metric case is chosen as $\phi_{d,A}(x)= d(x,A)$, and the constants $c_{d,A}$ are chosen independently of $d$ and are $0$ unless $|A|=2^s$ for some integer $s$ \cite[Thm 10.1.2]{Deza97}.

We can imagine following a similar method for diversities. Given a diversity $(X,\delta)$ with $|X|=n$, we select mappings $\phi_{\delta,A} \colon X \rightarrow \mathbb{R}$, and constants $c_{\delta,A} \geq 0$. Then we define 
\[
\hat{\delta}(B) = \sum_{A \subseteq X} c_{\delta,A} \max_{b_1,b_2 \in B} | \phi_{\delta, A}(b_1) - \phi_{\delta,A}(b_2) |,
\]
which is an $\ell_1$ diversity. The goal is then to define $\phi_{\delta,A}$ and $c_{\delta,A}$ to make $\hat{\delta}$ close to $\delta$.
We now describe some of the difficulties with making this work for diversities.

Here are three choices we could try for $\phi_{\delta,A}$.
\begin{enumerate}
\item $\phi_{\delta,A} = \delta_{A|\bar{A}}$ the split diversity for the split $A|\bar{A}$ \cite{bryant2014}.
\item $\phi_{\delta,A} = d(A,x) = \min_{a \in A} \delta(\{x,a\})$.
\item $\phi_{\delta,A} = \delta(A \cup \{x\})$.
\end{enumerate}
We consider each in turn.

(1). Any $\ell_1$ diversity can be written as a non-negative linear combination of the split diversities, so this choice of $\phi_{\delta,A}$ is just a restatement of the problem.

(2). This choice of $\phi_{\delta,A}$ only depends on $\delta$ through the values $\delta(\{x,a\}), x, a \in X$. In other words, the $\phi_{\delta,A}$ only depend on $\delta$ through its induced metric. In general $c_{\delta,A}$ may also depend on $\delta$, and in this is, as far as we know, a possible approach to the problem. But if we restrict ourselves to constants independent of $\delta$, $c_A$ for $A \subseteq X$, have a stronger restriction.
With this choice, any two diversities with the same induced metric will yield the same embedded diversity $\hat{\delta}$. This puts a constraint on how good the worst case distortion between $\delta$ and $\hat{\delta}$ can be. To see this, let $\delta_\rho$ be the discrete diversity on $n$ points that assigns value $1$ to all sets with two or more elements. Let $\delta_c$ be given by $\delta_2(A) = |A|-1$ for $|A|\geq 2$. $\delta_\rho$ and $\delta_c$ have the same induced metric, and so $\hat{\delta}$ would be identical for them. But $\delta_\rho(X)=1$ and $\delta_c(X)= n-1$ so $\hat{\delta}(X)$ must be off at least one of them by a factor of $\sqrt{n-1}$ or greater. So the best worst case distortion for this method is at least $\mathcal{O}(\sqrt{n})$. We don't know if this is attainable.

(3). The choice $\phi_{\delta,A}(x) = \delta(A \cup \{x\})$ looks promising in that it potentially uses all the information about the diversity.
We get
\[
\hat{\delta} (B) = \sum_{A \subseteq X} c_A \max_{b_1,b_2 \in B} \delta(A \cup \{b_1\}) - \delta(A \cup \{b_2\}).
\]
To get an idea of how this works, we apply it to the diversity $\delta_\rho$ from the previous case. Most terms are $0$ and  we get
\begin{eqnarray*}
\hat{\delta}(B) &=& \sum_{a \in X} c_{\delta_\rho,\{a\}} \max_{b_1,b_2 \in B} \delta_\rho(\{a,b_1\}) - \delta_\rho( \{ a,b_2 \}) \\
& = & \sum_{a \in X} c_{\delta_\rho,\{a\}} \mathbf{1}_{a \in B, |B|\geq 2} \\
& = &  \mathbf{1}_{|B|\geq 2} \sum_{a \in B} c_{\delta_\rho,\{a\}} = \delta(B) \sum_{a \in B} c_{\delta_\rho,\{a\}}.
\end{eqnarray*}
Suppose we have
\[
\hat{\delta}(B) \leq \delta(B) \leq g(n) \hat{\delta}(B).
\]
Taking $B=X$ this gives:
\[
\sum_{a \in X} c_{\delta_\rho,\{a\}} \leq 1.
\]
Taking $B= \{a_1, a_2\}$, $a_1 \not = a_2$ gives
\[
1 \leq g(n) (c_{\delta_\rho,\{a_1\}} + c_{\delta_\rho,\{a_2\}})
\]
Summing this over $\lfloor n/2 \rfloor$ pairs of distinct points in $X$ gives
\[
\lfloor n/2 \rfloor \leq g(n) \sum_{a \in X} c_{\delta_\rho,\{a\}} \leq g(n).
\]
So $g(n)$ is at least $\mathcal{O}(n)$. So this choice of $\phi_{\delta,A}$ cannot give us better results than the bound we already have.

\subsection{Bartal's Method}

Bartal's method \cite{bartal96,bartal98}, and as later refined in \cite{fakchar04}, relies on embedding a metric space $(X,d)$ in a random tree metric. Since tree metrics are $\ell_1$, a random tree metric is also $\ell_1$ (being a positive linear combination of tree metrics). These proofs work by randomly and hierarchically building a tree corresponding to the set $X$ such that tree metric dominates $d$ and, for each pair of nodes $x,y$, the expected distance is within $\mathcal{O}(\log n)$ of $d(x,y)$. 

The problem with this approach for diversities can be seen easily. We again use the diversity $\delta_\rho$ from above, where $\delta_\rho(A)= 1$ whenever $|A| \geq 2$. Since $\delta_\rho$ is a diameter diversity, we know it can be approximated with an $\ell_1$ diversity with distortion $\mathcal{O}(\log^2 n)$. But the best we can do with a phylogenetic diversity (the equivalent of a tree diversity) is much worse.

\begin{prop} Let $(X,\delta_\rho)$ be the diversity with $|X|=n$ and $\delta_\rho(A)= 1$ for all $|A|\geq 2$. Let $\delta_\tau$ be a random phylogenetic diversity with 
\[
\delta_\rho(A) \leq \delta_t (A), \ \ \, \mbox{and } \mathbb{E} \delta_t(A) \leq g(n) \delta_\rho(A)
\]
for all $A \subseteq X$. Then $g(n) \geq \lfloor n/2 \rfloor$.
\end{prop}
\begin{pf}
Let $L_t$ be the length of the tree $\delta_t$. By Corollary 5.18 and Theorem 5.6.2 in \cite{semple2003phylogenetics} there are $m= \lfloor n/2 \rfloor$ pairs of vertices such that the paths joining each pair are disjoint. Denoting these pairs $(a_1,b_1), \ldots, (a_m,b_m)$, we have that $L_t \geq \sum_{i=1}^m d(a_i,b_i) = \sum_{i=1}^m \delta_t(\{a_i,b_i\}) \geq m$. On the other hand,  
\[
m \leq \mathbb{E} L_t = \mathbb{E} \delta_t(X) \leq g(n) \delta_\rho(X) = g(n).
\]
So $g(n) \geq \lfloor n/2 \rfloor$.
\end{pf}

\section{Other ideas}\label{sec:otherIdeas}

\begin{enumerate}
\item Every diversity $(X,\delta)$ with a given induced metric $(X,d)$ satisfies the inequalities
\[\diam(A) \leq \delta(A) \leq \delta_S(A).\]
We have shown that Bourgain's approach embeds  $\diam$ with small distortion, but applying the metric result directly does appear to give an embedding for $\delta_S$. On the other hand, Bartal's approach embeds $\delta_S$, but does not successfully embed all diameter diversities. Is some hybrid approach possible, perhaps one which applies Bartal's approach but switches to a Bourgain approach for some subproblems?
\item Expander graphs have proven to be a good source for lower bounds on the distortion required to embed metrics in $\ell_1$. Expander hypergraphs are far less understood, but perhaps analogous arguments could be used to obtain lower bounds on the distortion of diversities, particularly given the connections between diversity embedding and hypergraph sparsest cut established in \cite{bryant2014}.
\item The embedding used to prove Theorem~\ref{PetersEmbedding} was derived completely from the induced metric. We know that we cannot obtain a distortion bound in this way that is better than $\Omega(\sqrt{n})$, as we showed in Part 2 of our discussion of Bourgain's approach. 
We conjecture that there exists an embedding based solely on the induced metric which nevertheless achieves an $\mathcal{O}(\sqrt{n})$ distortion bound for all diversities.
\end{enumerate}

\bibliographystyle{plain}

\end{document}